\documentclass[aos,MSNbibl,nameyear,dvips]{arximspdf}
\usepackage{graphicx}
\doi{10.1214/13-AOS1175E} 
\referstodoi{10.1214/13-AOS1175}
\volume{42}
\issue{2}
\pubyear{2014}
\firstpage{501}
\lastpage{508}

\makeatletter

\newproclaim{remark}[theorem]{Remark}
\newcommand{\argmin}{\mathop{\arg\min}}
\makeatother

\begin{document}
\begin{frontmatter}

\title{Discussion: ``A significance test for the lasso''}
\runtitle{Discussion}

\begin{aug}
\author{\fnms{Larry} \snm{Wasserman}\corref{}\ead[label=e1]{larry@stat.cmu.edu}}
\runauthor{L. Wasserman}
\affiliation{Carnegie Mellon University}
\address{Department of Statistics\\
Carnegie Mellon University\\
5000 Forbes Ave.\\
Pittsburgh, Pennsylvania 15213\\
USA\\
\printead{e1}} 
\pdftitle{Discussion of ``A significance test for the lasso''}
\end{aug}

\received{\smonth{12} \syear{2013}}



\end{frontmatter}

The paper by
Lockhart, Taylor, Tibshirani and Tibshirani (LTTT)
is an important advancement
in our understanding of inference
for high-dimension\-al regression.
The paper is a
tour de force,
bringing together
an impressive array of results,
culminating in a set of
very satisfying convergence results.
The fact that the test statistic
automatically balances
the effect of shrinkage and the effect
of adaptive variable selection is remarkable.

The authors make very strong assumptions.
This is quite reasonable: to make significant theoretical
advances in our understanding of complex procedures,
one has to begin with strong assumptions.
The following question then arises:
what can we do without these assumptions?

\section{The assumptions}

The assumptions in this paper---and in most theoretical papers on high-dimensional regression---have several components.
These include:

\begin{longlist}[(5)]
\item[(1)] The linear model is correct.
\item[(2)] The variance is constant.
\item[(3)] The errors have a Normal distribution.
\item[(4)] The parameter vector is sparse.
\item[(5)] The design matrix has very weak collinearity.
This is usually stated in the form of
incoherence, eigenvalue restrictions or incompatibility assumptions.
\end{longlist}

To the best of my knowledge, these assumptions are not testable when $p
> n$.
They are certainly a good starting place for
theoretical investigations but they are
indeed very strong.
The regression function
$m(x) = \mathbb{E}(Y|X=x)$ can be any function.
There is no reason to think it will be close to linear.
Design assumptions are also highly suspect.
High collinearity is the rule rather than the exception especially
in high-dimensional problems.
An exception is signal processing,
in particular compressed sensing,
where the user gets to construct the design matrix.
In this case, if the design matrix is filled with
independent random Normals, the design matrix will
be incoherent with high probability.
But this is a rather special situation.

None of this is meant as a criticism of the paper.
Rather, I am trying to motivate interest in the question
I asked earlier, namely:
what can we do without these assumptions?

%
\begin{remark}
It is also worth mentioning that even in
low-dimensional models and even if the model is correct,
model selection raises troubling issues that
we all tend to ignore.
In particular, variable selection
makes the minimax risk explode
[\citeauthor{LeePot08} (\citeyear{LeePot05,LeePot08})].
This is not some sort of pathological risk explosion, rather,
the risk is large in a neighborhood of 0,
which is a part of the parameter space we care about.
\end{remark}

\section{The assumption-free lasso}

To begin it is worth pointing out that the lasso
has a very nice assumption-free interpretation.

Suppose we observe
$(X_1,Y_1),\ldots, (X_n,Y_n)\sim P$
where
$Y_i\in\mathbb{R}$ and
$X_i \in\mathbb{R}^d$.
The regression function
$m(x) = \mathbb{E}(Y|X=x)$
is some unknown, arbitrary function.
We have no hope to estimate $m(x)$
nor do we have licence to impose assumptions on $m$.

But there is sound theory to justify the lasso
that makes virtually no assumptions.
In particular, I refer to
\citet{GreRit04} and
\citet{JudNem00}.

Let
${\mathcal L} = \{ x^t\beta\dvtx  \beta\in\mathbb{R}^d\}$
be the set of linear predictors.
For a given $\beta$,
define the predictive risk
\[
R(\beta) = E\bigl( Y - \beta^T X\bigr)^2,
\]
where $(X,Y)$ is a new pair.
Let us define the best, sparse, linear predictor
$\ell_*(x) = \beta_*^T x$
(in the $\ell_1$ sense)
where
$\beta_*$ minimizes
$R(\beta)$ over the set
$B(L) = \{ \beta\dvtx  \|\beta\|_1 \leq L\}$.
The lasso estimator $\hat\beta$
minimizes the empirical risk
$\hat R(\beta) = \frac{1}{n}\sum_{i=1}^n (Y_i - \beta^T X_i)^2$
over $B(L)$.
For simplicity, I will assume that
all the variables are bounded by $C$
(but this is not really needed).
We make no other assumptions:
no linearity,
no design assumptions
and no models.
It is now easily shown that
\[
R(\hat\beta) \leq R(\beta_*) + \sqrt{ \frac{8 C^2 L^4}{n} \log\biggl(
\frac{2p^2}{\delta} \biggr)}
\]
except on a set of probability at most $\delta$.

This shows that the predictive risk of the lasso
comes close to the risk of the best sparse linear predictor.
In my opinion, this explains why the lasso ``works.''
The lasso gives us a predictor with
a desirable property---sparsity---while being computationally tractable
and it comes close to the risk of the
best sparse linear predictor.

\section{Interlude: Weak versus strong modeling}

When developing new meth\-odology, I think it is useful
to consider three different stages of development:
\begin{longlist}[(3)]
\item[(1)] Constructing the method.
\item[(2)] Interpreting the output of the method.
\item[(3)] Studying the properties of the method.
\end{longlist}
I also think it is useful to
distinguish two types of modeling.
In \textit{strong modeling}, the model
is assumed to be true in all three stages.
In \textit{weak modeling}, the model
is assumed to be true for stage 1 but not for
stages 2 and 3.
In other words, one can use a model
to help construct a method.
But one does not have to assume the model is true
when it comes to interpretation or when
studying the theoretical properties of the method.
My discussion is guided by my preference
for weak modeling.

\section{Assumption-free inference: The HARNESS}

Here, I would like to discuss an approach I have been developing
with Ryan Tibshirani.
We call this:
{H}igh-dimensional Agnostic
{R}egression
{N}ot
{E}mploying
{S}tructure or
{S}parsity,
or, the HARNESS.
The method is a variant of the idea proposed in \citet{WasRoe09}.

The idea is to split the data into two halves.
${\mathcal D}_1$ and ${\mathcal D}_2$.
For simplicity, assume that $n$ is even so that each half has size
$m=n/2$.
From the first half ${\mathcal D}_1$, we select a subset of variables $S$.
The method is agnostic about how the variable selection is done.
It could be forward stepwise, lasso, elastic net or anything else.
The\vspace*{2pt} output of the first part of the analysis is
the subset of predictors
$S$ and an estimator
$\hat\beta= (\hat\beta_j\dvtx  j\in S)$.
The second half of the data ${\mathcal D}_2$ is used to provide
distribution-free inferences for the following questions:
\begin{longlist}[(3)]
\item[(1)] What is the predictive risk of $\hat\beta$?
\item[(2)] How much does each variable in $S$ contribute to the predictive risk?
\item[(3)] What is the best linear predictor using the variables in $S$?
\end{longlist}
All the inferences from ${\mathcal D}_2$
are interpreted as being conditional on ${\mathcal D}_1$.
(A~variation is to
use ${\mathcal D}_1$ only to produce $S$
and then construct the coefficients of the predictor from
${\mathcal D}_2$.
For the purposes of this discussion, we use $\hat\beta$
from ${\mathcal D}_1$.)

In more detail,
let
\[
R = \mathbb{E}\bigl| Y-X^T \hat\beta\bigr|,
\]
where the randomness is over the new pair $(X,Y)$;
we are conditioning on ${\mathcal D}_1$.
Note that in this section I have changed
the definition of $R$ to be on the absolute scale
which is more interpretable.
In the above equation, it is understood
that \mbox{$\hat\beta_j=0$} when $j\notin S$.
The first question refers to producing a estimate and confidence
interval for $R$
(conditional on ${\mathcal D}_1$).
The second question refers to inferring
\[
R_j = \mathbb{E}\bigl| Y-X^T \hat\beta_{(j)}\bigr| -
\mathbb{E}\bigl| Y-X^T \hat\beta\bigr|\vadjust{\goodbreak}
\]
for each $j\in S$,
where
$\beta_{(j)}$ is equal to $\hat\beta$
except that $\hat\beta_j$ is set to 0.
Thus,
$R_j$ is the risk inflation by excluding $X_j$.
The third question refers to
inferring
\[
\beta^* = \argmin_{\beta\in\mathbb{R}^k} \mathbb{E}\bigl(Y - X_S^T
\beta\bigr)^2
\]
the coefficient of the best linear predictor for the chosen model.
We call $\beta^*$ the
\textit{projected parameter}.
Hence, $x^T\beta^*$ is
the best linear approximation to $m(x)$
on the linear space spanned by the selected variables.

A consistent estimate of $R$ is
\[
\hat{R} = \frac{1}{m}\sum_{i=1}^m
\delta_i,
\]
where the sum is over ${\mathcal D}_2$, and
$\delta_i = |Y_i - X_i^T \hat\beta|$.
An approximate $1-\alpha$ confidence interval for $R$ is
$\hat R \pm z_{\alpha/2} s/\sqrt{m}$
where
$s$ is the standard deviation of the $\delta_i$'s.

The validity of this confidence interval is essentially
distribution-free.
In fact, if we want to be purely distribution-free
and avoid asymptotics,
we could instead define $R$ to be the median of
the law of $|Y-X^T \hat\beta|$.
Then the order statistics of the $\delta_i$'s can be used
in the usual way to get a finite sample, distribution-free
confidence interval for $R$.

Estimates and confidence intervals for $R_j$ can be obtained from
$e_1,\ldots, e_m$
where
\[
e_i = \bigl|Y_i - X^T \hat\beta_{(j)}\bigr|
- \bigl|Y_i - X^T \hat\beta\bigr|.
\]
Estimates and confidence intervals for $\beta_*$
can be obtained by standard least squares procedures
based on ${\mathcal D}_2$.
The steps are summarized in Figure~\ref{figharness}.

%
\begin{figure}[b]
\hrulefill\vspace*{6pt} 
\begin{center}
{\sf The HARNESS}\vspace*{6pt}
\end{center}
\flushleft{Input: data ${\mathcal D} = \{(X_1,Y_1),\ldots, (X_n,Y_n)\}$.}

\begin{longlist}[(3)]
\item[(1)] Randomly split the data into two halves
${\mathcal D}_1$ and ${\mathcal D}_2$.
\item[(2)] Use ${\mathcal D}_1$ to select a subset of variables $S$.
This can be forward stepwise, the lasso, or any other method.
\item[(3)] Let $R = \mathbb{E}( (Y-X^T \hat\beta)^2 | {\mathcal D}_1)$
be the predictive risk of the selected model
on a future pair $(X,Y)$, conditional on ${\mathcal D}_1$.
\item[(4)] Using ${\mathcal D}_2$ construct point estimates and confidence intervals
for $R$, $(R_j\dvtx  \hat\beta_j\neq0)$
and $\beta_*$.
\end{longlist}\vspace*{-4pt}
\hrulefill 
\caption{The steps in the HARNESS algorithm.}\label{figharness}
\end{figure}

The HARNESS bears some similarity to
POSI [\citet{Beretal13}]
which is another inference method for model selection.
They both eschew any assumption that the linear model
is correct.
But POSI attempts to make inferences that are valid over all
possible selected models while the HARNESS restricts
attention to the selected model.
Also, the HARNESS emphasizes predictive inferential statement.

\begin{figure}

\includegraphics{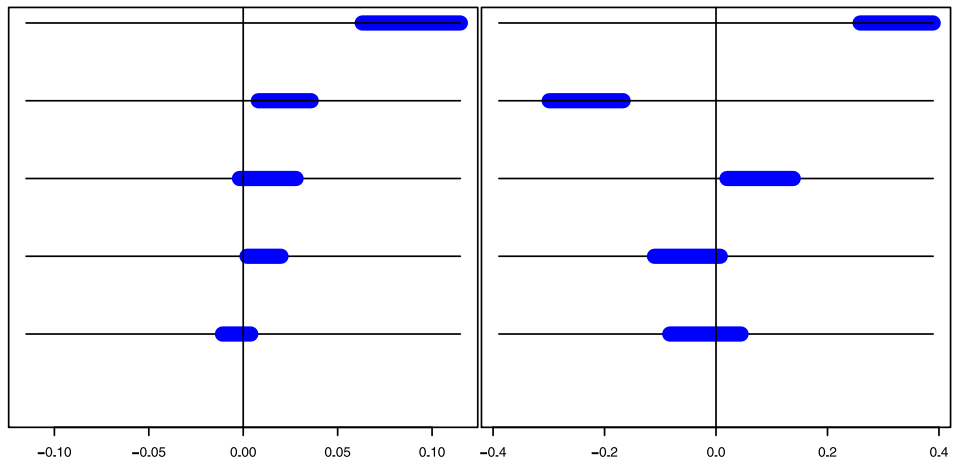}

\caption{Left plot: confidence intervals for $R_j$.
Right plot: confidence intervals for projected parameters.
From top down the variables are
Alcohol, Volatile\_Acidity, Sulphates,
Total\_Sulfur\_Dioxide and pH.}
\label{figwine}
\end{figure}

Here is an example using the wine dataset.
(Thanks to the authors for providing the data.)
Using the first half of the data,
we applied forward stepwise \mbox{selection} and used
$C_p$ to select a model.
The selected variables are
Alcohol, Volatile\_Acidity, Sulphates,
Total\_Sulfur\_Dioxide and pH.
A 95 percent confidence interval
for the predictive risk of the null model is
(0.65,0.70).
For the selected model, the confidence
interval for $R$ is
(0.46,0.53).
The (Bonferroni-corrected) 95 percent confidence intervals for the $R_j$'s
are shown in the first plot of Figure~\ref{figwine}.
The (Bonferroni-corrected) 95 percent confidence intervals for the
parameters of the projected model
are shown in the second plot in Figure~\ref{figwine}.

%

\section{The value of data-splitting}

Some statisticians are uncomfortable with data-splitting.
There are two common objections.
The first is that the inferences are random: if we repeat the procedure
we will get
different answers.
The second is that it is wasteful.

The first objection can be dealt with by
doing many splits and combining the information appropriately.
This can be done but is somewhat involved and will be described elsewhere.
The second objection is, in my view, incorrect.
The value of data splitting is that leads to
simple, assumption-free inference.
There is nothing wasteful about this.
Both halves of the data are being put to use.
Admittedly, the splitting leads to a loss of power
compared to ordinary methods
\textit{if the model were correct}.
But this is a false comparison since we are trying
to get inferences
without assuming the model is correct.
It is a bit like saying that nonparametric function
estimators have slower rates of convergence than parametric estimators.
But that is only because the parametric estimators invoke stronger assumptions.

\section{Conformal prediction}

Since I am focusing my discussion on
regression methods that make weak assumptions,
I would also like to briefly mention
Vladimir Vovk's theory of
conformal inference.
This is a completely distribution-free, finite sample
method for predictive regression.
The method is described in
\citet{VovGamSha05}
and
\citet{VovNouGam09}.
Unfortunately,
most statisticians seem to be unaware of this work
which is a shame.
The statistical properties (such as minimax properties)
of conformal prediction
were investigated in
\citet{Lei}, \citet{LeiWas}.

A full explanation of the method is beyond the scope of this
discussion but I do want to give the general idea
and say why it is related the current paper.
Given data
$(X_1,Y_1),\ldots,(X_n,Y_n)$,
suppose we observe a new $X$
and want to predict $Y$.
Let $y\in\mathbb{R}$ be an arbitrary real number.
Think of $y$ as a tentative guess at $Y$.
Form the augmented data set
\[
(X_1,Y_1),\ldots, (X_n,Y_n),
(X,y).
\]
Now we fit a linear model
to the augmented data
and compute residuals
$e_i$ for each of the $n+1$ observations.
Now we test
$H_0\dvtx  Y=y$.
Under $H_0$,
the residuals are invariant under permutations
and so
\[
p(y) = \frac{1}{n+1}\sum_{i=1}^{n+1} I\bigl(
|e_i| \geq|e_{n+1}|\bigr)
\]
is a distribution-free $p$-value for $H_0$.

Next, we invert the test:
let
$C = \{y\dvtx  p(y) \geq\alpha\}$.
It is easy to show that
\[
\mathbb{P}(Y\in C) \geq1-\alpha.
\]
Thus, $C$ is distribution-free, finite-sample
prediction interval for $Y$.
Like the \mbox{HARNESS},
the validity of the method does not depend
on the linear model being correct.
The set $C$ has the desired coverage probability no matter what
the true model is.
Both the HARNESS and conformal prediction use the linear model
as a device for generating predictions but neither
requires the linear model to be true
for the inferences to be valid.
[In fact, in the conformal approach, any method of generating residuals
can be used.
It does not have to be a linear model. See \citet{LeiWas}.]

One can also look at how the prediction interval $C$
changes as different variables are removed.\vadjust{\goodbreak}
This gives another assumption-free method to
explore the effects of predictors in regression.
Minimizing the length of the interval
over the lasso path can also be used as a
distribution-free method for choosing the regularization
parameter of the lasso.

On a related note,
we might also be interested in
assumption-free methods for the related task
of inferring graphical models.
For a modest attempt at this,
see \citet{WasKolRin}.

\section{Causation}

LTTT do not discuss causation.
But in any discussion of the assumptions
underlying regression,
causation is lurking just below the surface.
Indeed, there is a tendency to conflate
causation and inference.
To be clear: prediction, inference and causation
are three separate ideas.

Even if the linear model is correct,
we have to be careful how we interpret the parameters.
Many articles and textbooks describe $\beta_j$
as the change in $Y$ if $X_j$ is changed, holding the other covariates fixed.
This is incorrect.
In fact, $\beta_j$ is the \textit{change in our prediction of} $Y$
if $X_j$ is changed.
This may seem like nit-picking but this is the
very difference between association and causation.\looseness=-1

Causation refers to the change in $Y$
as $X_j$ is changed.
Association (prediction) refers to the change
\textit{in our prediction of} $Y$
as $X_j$ is changed.
Put another way,
prediction is about $\mathbb{E}(Y| \mbox{ observe } X=x)$ while
causation is about $\mathbb{E}(Y| \mbox{ set } X=x)$.
If $X$ is randomly assigned, they are the same.
Otherwise, they are different.
To make the causal claim, we have to include
in the model, every possible confounding variable
that could affect both $Y$ and $X$.
This complete causal model has the form
\[
Y = g(X,Z) + \varepsilon,
\]
where $Z=(Z_1,\ldots,Z_k)$
represents
all confounding variables in the world.
The relationship between $Y$ and $X$
alone is described as
\[
Y = f(X) + \varepsilon'.
\]
The causal effect---the change in $Y$
as $X_j$ is changed---is given by $\partial g(x,z)/\partial x_j$.
The association (prediction)---the change in our prediction of $Y$
as $X_j$ is changed---is given by $\partial f(x)/\partial x_j$.
If there are any omitted confounding variables,
then these will be different.

Which brings me back to the paper.
Even if the linear model
is correct,
we still have to exercise great caution in interpreting
the coefficients.
Most users of our methods
are nonstatisticians and are likely to interpret
$\beta_j$ causally no matter how many warnings we give.

\section{Conclusion}

LTTT have produced
a fascinating paper
that significantly advances
our understanding of
high-dimensional regression.
I expect there will be a flurry
of new research inspired by this paper.

My discussion has focused on the role of
assumptions.
In low-dimensional models,
it is relatively easy to
create methods that make few
assumptions.
In high-dimensional models,
low assumption inference is much more challenging.

I hope I have convinced the authors that
the low assumption world is worth exploring.
In the meantime, I congratulate the authors
on an important and stimulating paper.

\section*{Acknowledgments}
Thanks to Rob Kass, Rob Tibshirani and Ryan Tibshirani
for helpful comments.


%

\printaddresses

\end{document}